\documentclass{amsart}

\theoremstyle{definition}

\theoremstyle{remark}

\newtheorem{proposition2}{Proposition}
\newtheorem{theorem2}{Theorem}

\numberwithin{equation}{section}

\begin{document}

\title{An Isoperimetric Result on High-Dimensional Spheres}


\author{Leighton Pate Barnes}
\address{Department of Electrical Engineering, Stanford University, Stanford, CA 94305}
\email{lpb@stanford.edu}
\thanks{The work was supported in part by NSF award CCF-1704624 and by the Center for Science of Information (CSoI), an NSF Science and Technology Center, under grant agreement CCF-0939370.}

\author{Ayfer \"Ozg\"ur}
\address{Department of Electrical Engineering, Stanford University, Stanford, CA 94305}
\email{aozgur@stanford.edu}

\author{Xiugang Wu}
\address{Department of Electrical Engineering, Stanford University, Stanford, CA 94305}
\email{xwu@udel.edu}
\curraddr{Department of Electrical and Computer Engineering,
University of Delaware, Newark, DE 19716}

\subjclass[2010]{Primary 60D05}

\date{November 20, 2018}

\keywords{isoperimetric inequality, concentration of measure, high-dimensional geometry}

\commby{}

\begin{abstract}
We consider an extremal problem for subsets of high-dimensional spheres that can be thought of as an extension of the classical isoperimetric problem on the sphere. Let $A$ be a subset of the $(m-1)$-dimensional sphere $\mathbb{S}^{m-1}$, and let $\mathbf{y}\in \mathbb{S}^{m-1}$ be a randomly chosen point on the sphere. What is the measure of the intersection of the $t$-neighborhood of the point $\mathbf{y}$ with the subset $A$? We show that with high probability this intersection is approximately as large as the intersection that would occur with high probability if $A$ were a spherical cap of the same measure. 
\end{abstract}

\maketitle

\section{Introduction}


 Let $\mathbb{S}^{m-1}\subseteq \mathbb{R}^m$ denote the $(m-1)$-sphere of radius $R$, i.e., $$\mathbb{S}^{m-1}=\left\{\mathbf z\in\mathbb{R}^m:\|\mathbf z\|=R\right\},$$ equipped with the rotation invariant Haar measure $\mu$ that is normalized such that $$\mu(\mathbb{S}^{m-1}) = \frac{2\pi^{\frac{m}{2}}}{\Gamma(\frac{m}{2})}R^{m-1} \; .$$
This normalization corresponds to the usual surface area. Let $\mathbb{P}(A)$ denote the probability of a set $A$ (that we will always assume to be measurable) with respect to the corresponding Haar probability measure that is normalized such that $\mathbb{P}(\mathbb{S}^{m-1})=1$.
A spherical cap of angle $\theta$ and pole $\mathbf z_0$ is defined as a ball on $\mathbb{S}^{m-1}$ using the geodesic metric $\angle( \mathbf z,\mathbf y)=\arccos(\langle\mathbf z/R,\mathbf y/R\rangle)$, i.e.,
$$
\text{Cap}(\mathbf z_0,\theta)=\left\{\mathbf z\in \mathbb{S}^{m-1}: \angle(\mathbf z_0,\mathbf z)\leq \theta \right\}.
$$
We say that a set $A\subseteq \mathbb{S}^{m-1}$ has effective angle $\theta$ if $\mu(A)=\mu(C)$ with $C=\text{Cap}(\mathbf z_0, \theta)$ for some  $\mathbf z_0 \in \mathbb{S}^{m-1}$.

The classical isoperimetric inequality on the sphere implies that among all sets on the sphere with a given measure, the spherical cap has the smallest boundary, or more generally the smallest neighborhood \cite{levy},\cite{isoperimetric}. This is formalized as follows:

\begin{proposition2}[Isoperimetric Inequality]\label{prop:isoperimetry}
For any arbitrary set $A\subseteq \mathbb{S}^{m-1}$ such that $\mu(A)=\mu(C)$, where $C=\text{Cap}(\mathbf z_0,\theta)\subseteq \mathbb{S}^{m-1}$ is a spherical cap, it holds that $$\mu(A_t)\geq \mu(C_t), \ \forall t\geq 0,$$ where $A_t$ is the $t$-neighborhood of $A$ defined as $$A_t=\left\{\mathbf z\in \mathbb{S}^{m-1}: \min_{\mathbf z'\in A}\angle(\mathbf z, \mathbf z')\leq t\right\},$$ and similarly $$C_t=\left\{\mathbf z\in \mathbb{S}^{m-1}: \min_{\mathbf z'\in C}\angle(\mathbf z, \mathbf z')\leq t\right\}=\text{Cap}(\mathbf z_0,\theta+t).$$
\end{proposition2}

Another basic geometric phenomenon, this time occurring only in high-dimensions, is concentration of measure. On a high-dimensional sphere this manifests as most of the measure of the sphere concentrating around any equator. This is captured by the following proposition \cite{matou}.

\begin{proposition2}[Measure Concentration] \label{P:measurecon}
Given any $\epsilon>0$, there exists an $M(\epsilon)$ such that for any $m\geq M(\epsilon)$ and any $\  \mathbf z \in \mathbb{S}^{m-1}$,
\begin{align}
\mathbb P \left( \angle(\mathbf z,\mathbf Y) \in [\pi/2-\epsilon,\pi/2+\epsilon] \right)\geq 1-\epsilon , \end{align}
where $\mathbf Y\in\mathbb{S}^{m-1}$ is distributed according to the Haar probability measure.
\end{proposition2}

The above measure concentration result combined with the isoperimetric inequality immediately yields the following result:
\begin{proposition2}[Blowing-up Lemma]\label{P:blowup}
Let $A\subseteq \mathbb{S}^{m-1}$ be an arbitrary set with effective angle $\theta$. Then for any $\epsilon>0$ and $m$ sufficiently large,
\begin{equation}\label{eq:blowup0}
\mathbb{P}(A_{\frac{\pi}{2}-\theta+\epsilon})\geq 1-\epsilon.
\end{equation}
\end{proposition2}


An almost equivalent way to state the blowing-up lemma from Proposition~\ref{P:blowup} is the following: let $A\subseteq \mathbb{S}^{m-1}$ be an arbitrary set with effective angle $\theta>0$, then for any $\epsilon>0$ and sufficiently large $m$,
\begin{equation}\label{eq:isoperimetry}
\mathbb{P}\left(\mu\left(A\cap \text{Cap}\left(\mathbf Y,\frac{\pi}{2}-\theta+\epsilon\right)\right)>0\right)> 1-\epsilon,
\end{equation}
{where $\mathbf Y$ is distributed according to the normalized Haar measure on $\mathbb{S}^{m-1}$.} In words, if we take a $\mathbf y$ uniformly at random on the sphere and draw a spherical cap of angle slightly larger than $\frac{\pi}{2}-\theta$ around it, this cap will intersect the set $A$ with high probability. This statement is almost equivalent to \eqref{eq:blowup0} since the $\mathbf y$'s for which the intersection is not empty lie in the $\frac{\pi}{2}-\theta+\epsilon$-neighborhood of $A$. Note that this statement would trivially follow from measure concentration on the sphere (Proposition~\ref{P:measurecon}) if $A$ were known to be a spherical cap, and it holds for any $A$ due to the isoperimetric inequality in Proposition~\ref{prop:isoperimetry}.

Our main result is the following generalization of \eqref{eq:isoperimetry}.

{
\begin{theorem2}\label{thm:stongisoperimetry}
Given any $\epsilon>0$, there exists an $M(\epsilon)$ such that for any $m > M(\epsilon)$ the following is true. Let $A\subseteq \mathbb{S}^{m-1}$ be any arbitrary set with effective angle $\theta>0$, and let $V=\mu(\text{Cap}(\mathbf z_0, \theta) \cap \text{Cap}(\mathbf y_0, \omega))$ where  $\mathbf z_0,\mathbf y_0 \in \mathbb{S}^{m-1}$ with $\angle(\mathbf z_0,\mathbf y_0)=\pi/2$
and $\theta+\omega>\pi/2$. Then
\begin{equation*}
\mathbb{P}\left(\mu(A\cap \text{Cap}(\mathbf Y,\omega+\epsilon))> (1-\epsilon)V\right)\geq 1-\epsilon,
\end{equation*}
where $\mathbf Y$ is a random vector on $\mathbb{S}^{m-1}$ distributed according to the normalized Haar measure.
\end{theorem2}
}

If $A$ itself is a cap then the statement of Theorem 1 is straightforward and follows from the fact that $\mathbf y$ will be concentrated around the equator at angle $\pi/2$ from the pole of $A$ (Proposition~\ref{P:measurecon}). Therefore, as $m$ gets large, the intersection of the two spherical caps will be given by $V$  for almost all $\mathbf y$'s. The statement however is much stronger than this and holds for any arbitrary set $A$, analogous to the isoperimetric result in \eqref{eq:isoperimetry}. It states that no matter what the set $A$ is, if we take a random point on the sphere and draw a cap of angle slightly larger than $\omega$ for $\omega>\pi/2-\theta$, then with high probability the intersection of the cap with the set $A$ would be at least as large as the intersection we would get if $A$ were a spherical cap. 


The authors first encountered this problem while working on a geometric framework for proving capacity bounds in information theory. A version of this result on a spherical shell with nonzero thickness was used in \cite{ITtrans} to resolve the Gaussian case of an open problem from \cite{open}, and some partial results for a discrete version on the Hamming sphere were used in \cite{Allerton2017}.


When viewed as a concentration of measure result, it might seem as though Theorem \ref{thm:stongisoperimetry} could be proved with standard methods such as the concentration of Lipschitz functions \cite{isoperimetric}. We could not find a way to make this work, since the Lipschitz constant for the function $\mathbf y \mapsto \mu(A\cap \text{Cap}(\mathbf y,\omega+\epsilon))$ can increase exponentially in the dimension $m$. Even if one is only interested in the exponential order of the intersection measure, $\log(\mu(A\cap \text{Cap}(\mathbf y,\omega+\epsilon)))$ is equal to $-\infty$ when the sets are disjoint, so log composed with this function is not even continuous. However, if instead of the ``hard'' intersection measure that can be viewed as the indicator function $1_A$ convolved with the kernel $1_{\text{Cap}(\mathbf y,\omega+\epsilon)}$, we are interested in the ``soft'' intersection measure that is $1_A$ convolved with some smooth kernel, then these techniques could possibly be applied. In particular, $\sqrt{\frac{1}{m}\log{\frac{1}{P_t1_A}}}$ exhibits the correct order of Lipschitz continuity where $P_t$ denotes the heat semigroup over time $t$ on the sphere. Once concentration of $\frac{1}{m}\log{P_t1_A}$ has been established, a rearrangement theorem such as Theorem 4.1 from \cite{continuous} could be used to establish that the smallest value it can concentrate around is given when $A$ is a spherical cap \footnote{Special thanks to Ramon van Handel for pointing out that concentration of Lipschitz functions can indeed be used in this ``soft'' intersection case.}.

\section{Rearrangement on the Sphere}

Our main tool for proving  Theroem~\ref{thm:stongisoperimetry} for arbitrary $A$ is the symmetric decreasing rearrangement of functions on the sphere, along with a version of the Riesz rearrangement inequality on the sphere due to Baernstein and Taylor \cite{baernstein}.

For any measurable function $f:\mathbb{S}^{m-1} \to \mathbb{R}$ and pole $\mathbf z_0$, the symmetric decreasing rearrangement of $f$ about $\mathbf z_0$ is defined to be the function $f^*:\mathbb{S}^{m-1} \to \mathbb{R}$ such that $f^*(\mathbf y)$ depends only on the angle $\angle(\mathbf y,\mathbf z_0)$, is nonincreasing in $\angle (\mathbf y,\mathbf z_0)$, and has super-level sets of the same size as $f$, i.e.

$$\mu\big(\{\mathbf y:f^*(\mathbf y)>d\}\big) = \mu\big(\{\mathbf y:f(\mathbf y)>d\}\big)$$
for all $d$. The function $f^*$ is unique except for on sets of measure zero.

One important special case is when the function $f=1_A$ is the characteristic function for a subset $A$. The function $1_A$ is just the function such that
$$1_A(\mathbf y) = \begin{cases} 1 & \mathbf y \in A \\ 0 &  \text{otherwise.} \end{cases} $$
In this case, $1_A^*$ is equal to the characteristic function associated with a spherical cap of the same size as $A$. In other words, if $A^*$ is a spherical cap about the pole $\mathbf z_0$ such that $\mu(A^*) = \mu(A)$, then $1_A^* = 1_{A^*}$.

\begin{theorem2}[Baernstein and Taylor \cite{baernstein}]\label{thm:riesz}
Let $K$ be a nondecreasing bounded measurable function on the interval $[-1,1]$. Then for all functions $f,g \in L^1(\mathbb{S}^{m-1})$,

\begin{eqnarray*}
& &\int_{\mathbb{S}^{m-1}} \left( \int_{\mathbb{S}^{m-1}} f(\mathbf z)K\left(\langle \mathbf z/R,\mathbf y /R\rangle\right) d\mathbf z \right) g(\mathbf y)d\mathbf{y} \\
& & \leq \int_{\mathbb{S}^{m-1}} \left( \int_{\mathbb{S}^{m-1}} f^*(\mathbf z)K\left(\langle \mathbf z/R, \mathbf y/R \rangle\right) d\mathbf z \right) g^*(\mathbf y)d\mathbf{y}.
\end{eqnarray*}
\end{theorem2}

{
\section{Proof of Theorem~\ref{thm:stongisoperimetry}}

In order to prove Theorem~\ref{thm:stongisoperimetry}, we will apply Theorem \ref{thm:riesz} by choosing $f,K$ such that the inner integral
$$ \int_{\mathbb{S}^{m-1}} f(\mathbf z)K\left(\langle \mathbf z/R,\mathbf y /R\rangle\right) d\mathbf{z} = \mu(A\cap \text{Cap}(\mathbf y,\omega+\epsilon)).$$
To do this, given an arbitrary set $A$, we set $f = 1_A$ and
$$K(\cos\alpha) = \begin{cases} 1 &  0\leq \alpha \leq \omega + \epsilon\\ 0 &  \omega + \epsilon < \alpha \leq \pi. \end{cases}$$
Note that $K$ is nondecreasing, bounded, and measurable. Furthermore, the product $f(\mathbf z)K\left(\langle \mathbf z/R,\mathbf y/R \rangle\right)$ is one precisely when $\mathbf z \in A$ and $\angle(\mathbf z,\mathbf y) \leq \omega + \epsilon$, and it is zero otherwise. Thus the integral
\begin{equation}\label{eq:inner}
\psi(\mathbf y) = \int_{\mathbb{S}^{m-1}} f(\mathbf z)K\left(\langle \mathbf z/R,\mathbf y/R \rangle\right) d\mathbf{z}
\end{equation}
is exactly the measure of the set $A \cap \text{Cap}(\mathbf y,\omega + \epsilon)$. 

We will use test functions $g$ that are also characteristic functions. Let $g=1_C$ where $C = \{\mathbf y:\psi(\mathbf y) > d\}$ for some $d$ (i.e. $C$ is a super-level set). For a fixed measure $\mu(C)$, the left-hand side of the inequality from Theorem \ref{thm:riesz} will be maximized by this choice of $C$. With this choice we have the following equality:
\begin{eqnarray*}
\int_{\mathbb{S}^{m-1}} \psi(\mathbf y) 1_C(\mathbf y)d\mathbf{y} & = & \int_{\mathbb{S}^{m-1}} \psi^*(\mathbf y) 1^*_C(\mathbf y)d\mathbf{y} \\
& = & \int_{C^*} \psi^*(\mathbf{y})d\mathbf {y} .
\end{eqnarray*}
This follows from the layer-cake decomposition for any non-negative and measurable function $\psi$ in that
\begin{eqnarray}
\int_{\mathbb{S}^{m-1}} \psi(\mathbf{y}) 1_C(\mathbf{y})d\mathbf{y} & = & \int_C \psi(\mathbf{y})d\mathbf{y} \nonumber \\
& = &\int_C \int_0^\infty 1_{\{\psi(\mathbf{y})>t\}}dt d\mathbf{y} \nonumber\\
& = &\int_0^\infty \int_C 1_{\{\psi(\mathbf y)>t\}}d\mathbf{y} dt \nonumber\\
& = &\int_0^d \int_{C} 1_{\{\psi(\mathbf y)>t\}}d\mathbf{y} dt   + \int_d^\infty \int_{C} 1_{\{\psi(\mathbf y)>t\}}d\mathbf{y} dt \nonumber\\
& = & \int_0^d \int_{C^*} 1_{\{\psi^*(\mathbf y)>t\}}d\mathbf{y} dt   + \int_d^\infty \int_{C^*} 1_{\{\psi^*(\mathbf y)>t\}}d\mathbf{y} dt \nonumber\\
& = & \int_{C^*} \psi^*(\mathbf{y})d\mathbf{y} \; . \label{eq:layer_cake}
\end{eqnarray}
Using this equality and our choices for $f,g,K$ above we will rewrite the inequality from Theorem \ref{thm:riesz} as
\begin{equation}\label{eq:ineq}
\int_{C^*} \psi^*(\mathbf y)d\mathbf{y} \leq \int_{C^*} \bar{\psi}(\mathbf y) d\mathbf{y}
\end{equation}
where
$$\bar{\psi}(\mathbf y) = \int_{\mathbb{S}^{m-1}} f^*(\mathbf z)K\left(\langle \mathbf z/R,\mathbf y/R \rangle\right) d\mathbf{z} \; .$$
Note that $\bar\psi(\mathbf y)$ is exactly $\mu(A^* \cap \text{Cap}(\mathbf y,\omega + \epsilon))  .$

Note that both $\psi^*(\mathbf y)$ and $\bar{\psi}(\mathbf y)$ are spherically symmetric. More concretely, they both depend only on the angle $\angle(\mathbf y,\mathbf z_0)$, so in an abuse of notation we will write $\bar{\psi}(\alpha)$ and $\psi^*(\alpha)$ where $\alpha = \angle(\mathbf y,\mathbf z_0)$.

For convenience we will define a measure $\nu$ by
$$d\nu(\phi) = A_{m-2}(R\sin\phi)R d\phi$$
where $A_m(R)$ denotes the Haar measure of the $m$-sphere with radius $R$. We do this so that an integral like
$$\int_{\mathbb{S}^{m-1}} \psi^* d\mathbf{y} = \int_0^\pi \psi^*(\phi)A_{m-2}(R\sin\phi)Rd\phi$$ can be expressed as
$$\int_0^\pi \psi^* d\nu \; .$$

We are now ready to prove Theroem~\ref{thm:stongisoperimetry}. Proposition \ref{P:measurecon} implies that for any $0 < \epsilon < 1$, there exists an $M(\epsilon)$ such that for $m > M(\epsilon)$ we have
\begin{eqnarray} \label{eq:con1}
\mathbb P \left( \angle(\mathbf z_0,\mathbf Y) \in [\pi/2-\epsilon,\pi/2+\epsilon] \right) & \geq & 1-\frac{\epsilon^2}{2} \; .
\end{eqnarray}
The constant $M(\epsilon)$ is determined only by the concentration of measure phenomenon cited above, and it does not depend on any parameters in the problem other than $\epsilon$. From now on, let us restrict our attention to dimensions $m > M(\epsilon)$. Due to the triangle inequality for the geodesic metric, for $\mathbf y$ such that $\angle(\mathbf z_0,\mathbf y) \in [\pi/2-\epsilon,\pi/2+\epsilon]$ we have
$$
A^*\cap \text{Cap}(\mathbf y_0, \omega) \subseteq A^* \cap \text{Cap}(\mathbf y,\omega+\epsilon)
$$
where $\mathbf y_0$ is such that $\angle(\mathbf z_0,\mathbf y_0)=\pi/2$. Therefore,
\begin{equation} \label{eq:greaterV}
 \bar{\psi}(\angle(\mathbf z_0,\mathbf y)) = \mu(A^* \cap \text{Cap}(\mathbf y, \omega + \epsilon)) \geq V
\end{equation}
for all for $\mathbf y$ such that $\angle(\mathbf z_0,\mathbf y) \in [\pi/2-\epsilon,\pi/2+\epsilon]$ and

\begin{eqnarray}
\mathbb P\left(  \bar{\psi}(\mathbf{Y})\geq V\right) & = &
\mathbb P\left( \mu(A^*\cap \text{Cap}(\mathbf Y,\omega+\epsilon))\geq V \right) \nonumber \\
& \geq & 1-\frac{\epsilon^2}{2} \nonumber\\
& \geq & 1 - \frac{\epsilon}{2} \label{eq:con2} \; .
\end{eqnarray}

To prove the lemma, we need to show that
\begin{equation} \label{eq:lem31goal1}
\mathbb P\left(  \psi(\mathbf Y) > (1-\epsilon)V \right) = \mathbb P\left( \mu(A\cap \text{Cap}(\mathbf Y,\omega+\epsilon)) > (1-\epsilon)V \right) \geq 1-\epsilon
\end{equation}
for any arbitrary set $A\subset \mathbb{S}^{m-1}$. Recall that by the definition of a decreasing symmetric rearrangement, we have
\begin{equation*}
\mathbb P\left(  \psi^*(\mathbf Y)>d \right) = \mathbb P\left( \psi(\mathbf Y)>d \right)
\end{equation*}
for any threshold $d$ and this implies
\begin{equation} \label{eq:levelseteq}
\mathbb P\left(  \psi^*(\mathbf Y) \leq (1-\epsilon)V \right) = \mathbb P\left( \psi(\mathbf Y) \leq (1-\epsilon)V \right) \; .
\end{equation}
 Therefore, the desired statement in \eqref{eq:lem31goal1} can be equivalently written as
\begin{equation} \label{eq:lem31step1}
\mathbb P\left(  \psi^*(\mathbf Y) \leq (1-\epsilon)V \right) \leq \epsilon.
\end{equation}

{Turning to proving \eqref{eq:lem31step1}, recall that by the definition  of a decreasing symmetric rearrangement, $\psi^*(\alpha)$ is nonincreasing over the interval $0\leq \alpha \leq \pi$. Let $\beta$ be the smallest value such that $\psi^*(\beta) = (1-\epsilon)V$, or more explicitly, $$\beta = \inf\{\alpha | \psi^*(\alpha) \leq (1-\epsilon)V\} \; .$$ If  $\beta \geq \pi/2 + \epsilon$, then \eqref{eq:lem31step1} would follow trivially from \eqref{eq:con1} and the fact that $\psi^*(\alpha)$ would be greater than $(1-\epsilon)V$ for all $0<\alpha<\pi/2 + \epsilon$. We will therefore assume that $0< \beta < \pi/2 + \epsilon$. It remains to show that even if this is the case, we have \eqref{eq:lem31step1}.}

By the definition of $\beta$ and the fact that $\psi^*$ is nonincreasing,
\begin{eqnarray}
\mathbb P\left(  \psi^*(\mathbf Y) \leq (1-\epsilon)V \right) & = & \frac{1}{A_{m-1}(R)}\int^\pi_\beta d\nu \nonumber\\
& = & \frac{1}{A_{m-1}(R)}\int^{\max\{\beta,\frac{\pi}{2}-\epsilon\}}_\beta d\nu + \frac{1}{A_{m-1}(R)}\int^{\frac{\pi}{2}+\epsilon}_{\max\{\beta,\frac{\pi}{2}-\epsilon\}} d\nu \nonumber \\
& + & \frac{1}{A_{m-1}(R)}\int^\pi_{\frac{\pi}{2}+\epsilon} d\nu \label{eq:psilessthanV2} \; .
\end{eqnarray}
To bound the first and third terms of \eqref{eq:psilessthanV2} note that
\begin{eqnarray}
 \frac{1}{A_{m-1}(R)}\int^{\max\{\beta,\frac{\pi}{2}-\epsilon\}}_\beta d\nu + \frac{1}{A_{m-1}(R)}\int^\pi_{\frac{\pi}{2}+\epsilon} d\nu & \leq & \frac{\epsilon^2}{2} \label{eq:conineq0}\\
& \leq & \frac{\epsilon}{2} \label{eq:conineq}
\end{eqnarray}
as a consequence of \eqref{eq:con1}. In order to bound the second term, we establish the following chain of (in)equalities which will be justified below.
\begin{align}
\frac{1}{A_{m-1}(R)}\int^\pi_{\frac{\pi}{2}+\epsilon}d\nu & \geq \frac{1}{(1-\epsilon)V A_{m-1}(R)}\int^\pi_{\frac{\pi}{2}+\epsilon}(\psi^*-\bar{\psi})d\nu \label{eq:step1}\\
& =  \frac{1}{(1-\epsilon)V A_{m-1}(R)}\int_0^{\frac{\pi}{2}+\epsilon}(\bar{\psi}-\psi^*)d\nu \label{eq:step2}\\
& \geq  \frac{1}{(1-\epsilon)V A_{m-1}(R)}\int_{\beta}^{\frac{\pi}{2}+\epsilon}(\bar{\psi}-\psi^*)d\nu \label{eq:step3}\\
& \geq  \frac{\epsilon}{(1-\epsilon)A_{m-1}(R)}\int_{\max\{\beta,\frac{\pi}{2}-\epsilon\}}^{\frac{\pi}{2}+\epsilon}d\nu \label{eq:step4} \\
& \geq \frac{\epsilon}{A_{m-1}(R)}\int_{\max\{\beta,\frac{\pi}{2}-\epsilon\}}^{\frac{\pi}{2}+\epsilon}d\nu \label{eq:step5}
\end{align}
Combining \eqref{eq:step5} with \eqref{eq:conineq0} reveals that the second term in \eqref{eq:psilessthanV2} is also bounded by $\epsilon/2$, therefore
$$\mathbb P\left(  \psi^*(\mathbf Y) \leq (1-\epsilon)V \right)$$
must be bounded by $\epsilon$ which implies the lemma.

The first inequality \eqref{eq:step1} is a consequence of the fact that over the range of the integral, $\psi^*$ is less than or equal to $(1-\epsilon)V$ and $\bar{\psi}$ is non-negative. The equality in \eqref{eq:step2} follows from
$$\int_0^\pi \psi^*d\nu = \int_0^\pi \bar{\psi}d\nu \; ,$$
which is itself a consequence of \eqref{eq:layer_cake} with $C = \mathbb{S}^{m-1}$ and
\begin{align}
\int_{\mathbb{S}^{m-1}} \psi(\mathbf y) d\mathbf y = \int_{\mathbb{S}^{m-1}} \int_{\mathbb{S}^{m-1}} f(\mathbf z) & K(\langle\mathbf z/R, \mathbf y/R\rangle)d\mathbf zd\mathbf y \nonumber\\
& = \int \int K(\langle \mathbf y,\mathbf z\rangle)d\mathbf y f(\mathbf z)d\mathbf z \nonumber\\
& =  \int \mu(\text{Cap}(\mathbf y,\omega)) f(\mathbf z) d\mathbf z \nonumber\\
& =  \mu(\text{Cap}(\mathbf y,\omega))\mu(A) \nonumber\\
& =  \int \mu(\text{Cap}(\mathbf y,\omega)) f^*(\mathbf z) d\mathbf z\nonumber \\
& =  \int \int f^*(\mathbf z)K(\langle\mathbf z/R,\mathbf y/R \rangle) d\mathbf zd\mathbf y = \int_{\mathbb{S}^{m-1}} \bar{\psi}(\mathbf y) d\mathbf y \label{eq:totalmeas}
 \; .
\end{align}
Next we have \eqref{eq:step3} which is due to the rearrangement inequality \eqref{eq:ineq} when $C$ is the super-level set $\{{\mathbf y} \; | \; \psi(\mathbf y) > (1-\epsilon)V\}$. By the definition of a symmetric decreasing rearrangement, $\mu(\{{\mathbf y} \; | \; \psi(\mathbf y) > (1-\epsilon)V\}) = \mu(\{{\mathbf y} \; | \; \psi^*(\mathbf y) > (1-\epsilon)V\})$, and the set on the right-hand side is an open or closed spherical cap of angle $\beta$. Thus $C^*$ is a spherical cap with angle $\beta$ and the rearrangement inequality \eqref{eq:ineq} gives $$\int_0^\beta \psi^* d\nu \leq \int_0^\beta \bar\psi d\nu \; .$$
Finally, for the inequality \eqref{eq:step4}, we first replace the lower integral limit with $\max\{\beta,\pi/2-\epsilon\} \geq \beta$. Then $\bar{\psi} \geq V$ over the range of the integral due to \eqref{eq:greaterV}.  Additionally, $\psi^* \leq (1-\epsilon)V$ over the range of the integral, and the inequality follows.
}


\bibliographystyle{amsplain}

\end{document}